\documentclass[11pt,a4paper,reqno]{amsart}

\usepackage{amssymb,amsmath}
\usepackage{latexsym}
\usepackage{graphics}
\usepackage{euscript}
\usepackage{mathrsfs}
\usepackage[dvips]{curves}
\usepackage{tikz}
\usepackage{pgf}
\usetikzlibrary{arrows,decorations.pathmorphing,backgrounds,positioning,fit}
\usepackage{epsf}
\usepackage{xcolor}
\usepackage[percent]{overpic}
\usepackage{setspace}
\usepackage[colorlinks=true,linktocpage=true,pagebackref=false, urlcolor=black, citecolor=black,linkcolor=black]{hyperref}
\usepackage{multirow}

\def\sqbullet{\raise.2ex\hbox{\vrule width 3.5pt height 3.5pt}}

{\em}

\newcounter{substep}
\def\thesubstep{\arabic{substep}}

\newenvironment{substeps}[1]{%
\refstepcounter{substep}\noindent{(\ref{#1}.\thesubstep)\ }\ }%
{\em}

\newcounter{subsubstep}
\def\thesubsubstep{\arabic{subsubstep}}

{\em}

\newtheorem{mthm}{Theorem}
\newtheorem*{mainthm}{Main Theorem}
\newtheorem{mprop}[mthm]{Proposition}
\newtheorem{thm}{Theorem}[section]

\newtheorem{cor}[thm]{Corollary}
\newtheorem{lem}[thm]{Lemma}

\theoremstyle{definition}

\newtheorem{define}[thm]{Definition}

\newtheorem*{property}{Property}

\theoremstyle{remark}
\newtheorem{remark}[thm]{Remark}

 \newcommand{\R}{{\mathbb R}}

\newcommand{\sph}{{\mathbb S}} 
 \newcommand{\PP}{{\mathbb P}}

\newcommand{\pol}{{\EuScript K}}

\newcommand{\p}{{\EuScript P}}

\newcommand{\Ss}{{\EuScript S}}

\newcommand{\Aa}{{\EuScript A}}

\newcommand{\Bb}{{\EuScript B}}

\newcommand{\Ii}{{\EuScript I}}

\newcommand{\Ff}{{\EuScript F}}

\newcommand{\Cc}{{\EuScript C}}

\newcommand{\Ee}{{\EuScript E}}

\newcommand{\tildebaja}{{\raise.17ex\hbox{$\scriptstyle\sim$}}}

\newcommand{\Int}{\operatorname{Int}}
\newcommand{\im}{\operatorname{Im}}

\newcommand{\cl}{\operatorname{Cl}}

\newcommand{\id}{\operatorname{id}}

\newcommand{\conv}[2]{\vec{{\mathfrak C}}_{#2}(#1)}

\newcommand{\x}{{\tt x}}

\newcommand{\ol}{\overline}

\numberwithin{equation}{section}

\title[Complements of convex polyhedra as polynomial images of $\R^3$]{On the complements of 3-dimensional convex polyhedra as polynomial images of $\R^3$}

\author{Jos\'e F. Fernando}
\address{Departamento de \'Algebra, Facultad de CC. Matem\'aticas, Universidad Complutense de Madrid, 28040 MADRID (SPAIN)}
\curraddr{}
\email{josefer@mat.ucm.es}
\thanks{First author is supported by Spanish GR MTM2011-22435, while the second is a external collaborator.
}

\author{Carlos Ueno}
\address{Departamento de Matem\'aticas, IES La Vega de San Jos\'e, Paseo de San Jos\'e, s/n, Las Palmas de Gran Canaria, 35015 LAS PALMAS (SPAIN).
}
\email{cuenjac@gmail.com}
\begin{document}

\begin{abstract}
Let $\pol\subset\R^n$ be a convex polyhedron of dimension $n$. Denote $\Ss:=\R^n\setminus\pol$ and let $\ol{\Ss}$ be its closure. We prove that for $n=3$ the semialgebraic sets $\Ss$ and $\ol{\Ss}$ are polynomial images of $\R^3$. The former techniques cannot be extended in general to represent the semialgebraic sets $\Ss$ and $\ol{\Ss}$ as polynomial images of $\R^n$ if $n\geq4$. 
\end{abstract}
\date{07/03/2014}

\maketitle

\section*{Introduction}

A map $f:=(f_1,\ldots,f_m):\R^n\to\R^m$ is a \em polynomial map \em if its components $f_k\in\R[\x]:=\R[\x_1,\ldots,\x_n]$ are polynomials. Analogously, $f$ is a \em regular map \em if its components can be represented as quotients $f_k=g_k/h_k$ of two polynomials $g_k,h_k\in\R[\x]$ such that $h_k$ never vanishes on $\R^n$. During the last decade we have approached the following problem: 

\vspace*{2mm}
\centerline{\em To determine which subsets $\Ss\subset\R^m$ are polynomial or regular images of $\R^n$\em.} 

\vspace*{1mm}
\noindent We refer to \cite{g} for the first proposal of studying this problem and some particular related ones like the ``quadrant problem'' \cite{fg1}. By Tarski-Seidenberg's principle \cite[1.4]{bcr} the image of an either polynomial or regular map is a semialgebraic set. A subset $\Ss\subset\R^n$ is \em semialgebraic \em when it has a description by a finite boolean combination of polynomial equations and inequalities, which we will call a \em semialgebraic \em description.

The effective representation of a subset $\Ss\subset\R^m$ as polynomial or regular image of $\R^n$ reduces the study of certain classical problems in Real Geometry to its study in $\R^n$ with the advantage of avoiding contour conditions. Examples of these problems are Optimization or Positivstellens\"atze certificates \cite{fg2,fu2}. The latter can be generalized to non-necessarily closed basic semialgebraic sets (like the exterior of a convex polyhedron) using this type of representations. Polynomial representations are advantageous to regular representations because they avoid denominators and in this work we show that the quite natural polynomial constructions devised in \cite{fg2,u2} still work for dimension $3$ in the unbounded (remaining) case. For higher dimension, our results in \cite{fg2} concerning regular images are still the best for the unbounded case.

We are far from solving the representation problems stated above in its full generality but we have developed significant progresses in two ways:

\noindent{\em General conditions.} By obtaining general conditions that must satisfy a semialgebraic subset $\Ss\subset\R^m$ that is either a polynomial or a regular image of $\R^n$ (see \cite{fg2,fu,u1}). The most remarkable one states that the set of points at infinity of a polynomial image of $\R^n$ is connected.

\noindent{\em Ample families.} By showing constructively that ample families of significant semialgebraic sets are either polynomial or regular images of $\R^n$ (see \cite{f1,fg1,fgu1,fu2,u2}).

A distinguished family of semialgebraic sets is the one constituted by those whose boundary is piecewise linear, that is, semialgebraic sets that admit a semialgebraic description involving just linear equations. Many of them cannot be polynomial or regular images of $\R^n$, but it seems natural to wonder what happens with convex polyhedra, their interiors (as topological manifolds with boundary), their complements and the complements of their interiors. As the $1$-dimensional case is completely determined in \cite{f1} we assume dimension $\geq2$. 

We proved in \cite{fgu1} that all $n$-dimensional convex polyhedra and their interiors are regular images of $\R^n$. As many convex polyhedra are bounded and the images of nonconstant polynomial maps are unbounded, the suitable approach there was to consider regular maps. Concerning the representation of unbounded polygons as polynomial images see \cite{u3}. In \cite{fu2} we prove that the complement $\R^n\setminus\pol$ of a convex polyhedron $\pol\subset\R^n$ that does not disconnect $\R^n$ and the complement $\R^n\setminus\Int\pol$ of its interior are regular images of $\R^n$. If $\pol$ is moreover bounded or has dimension $d<n$, then $\R^n\setminus\pol$ and $\R^n\setminus\Int\pol$ are polynomial images of $\R^n$. The techniques we developed in \cite{fu2} can be squeezed out to represent as polynomial images of $\R^3$ the semialgebraic sets $\R^3\setminus\pol$ and $\R^3\setminus\Int\pol$ if $\pol\subset\R^3$ is a $3$-dimensional unbounded convex polyhedron. Analogous results appear in \cite{u2} for dimension $2$.

\begin{mainthm}\label{main}
Let $\pol\subset\R^3$ be a $3$-dimensional unbounded convex polyhedron that does not disconnect $\R^3$. Then $\R^3\setminus\pol$ and $\R^3\setminus\Int\pol$ are polynomial images of $\R^3$. 
\end{mainthm}

A convex polyhedron $\pol\subset\R^n$ disconnects $\R^n$ if and only if it is a \em layer \em, that is, it is affinely equivalent to $[-a,a]\times\R^{n-1}$ for some $a\in\R$, which reduces to a hyperplane if $a=0$. These layers are particular cases of \em degenerate \em convex polyhedra that do not have vertices. We will refer to convex polyhedra that have at least one vertex as \em non-degenerate\em.

\subsection*{Trimming positions with respect to a facet}
The proof of the Main Theorem is based on a ``placing problem'' that involves the following definitions introduced in \cite[4.2 \& 6.1]{fu2}; we simplify them here to ease the discussion. Consider the fibration of $\R^n$ induced by the projection $\pi_{n}:\R^n\to\R^{n-1}\times\{0\},\ (x_1,\dots,x_n)\mapsto (x_1,\dots,x_{n-1},0)$. The fiber $\pi_n^{-1}(a,0)$ is a parallel line to the vector $\vec{e}_n:=(0,\ldots,0,1)$ for each $a\in\R^{n-1}$. Given an $n$-dimensional convex polyhedron $\pol\subset\R^n$, the intersection $\Ii_a:=\pi_n^{-1}(a,0)\cap\pol$ can be either empty or a closed interval (either bounded or unbounded). Let $\pol\subset\R^n$ be an $n$-dimen\-sion\-al convex polyhedron and let $\Ff$ be one of its facets (faces of dimension $n-1$). Consider the set ${\mathfrak A}_\pol:=\{a\in\R^{n-1}:\ \Ii_a\neq\varnothing,\ (a,0)\not\in\Ii_a\}$. We say that:

(1) $\pol$ is in \em first trimming position with respect to the facet $\Ff$ \em if:
\begin{itemize}
\item[(i)] $\Ff\subset\{x_{n-1}=0\}$ and $\pol\subset\{x_{n-1}\leq0\}$.
\item[(ii)] For all $a\in\R^{n-1}$ the interval $\Ii_a$ is bounded. 
\item[(iii)] The set ${\mathfrak A}_\pol$ is bounded.
\end{itemize}

(2) $\pol$ is in \em second trimming position with respect to the facet $\Ff$ \em if:
\begin{itemize}
\item[(i)] $\Ff\subset\{x_n=0\}$ and $\pol\subset\{x_n\leq0\}$.
\item[(ii)] The set ${\mathfrak A}_\pol$ is bounded.
\end{itemize}

A strategy to place an unbounded convex polyhedron $\pol\subset\R^n$ in a trimming position with respect to a facet requires to characterize when the set ${\mathfrak A}_\pol$ is bounded and this is done in Lemma \ref{proysec}. Using it we get for the $3$-dimensional case the following placing results.

\begin{mprop}[$3$-dimensional polyhedra and first trimming position]\label{1trim} 
Let $\pol\subset\R^3$ be a non-degenerate $3$-dimensional unbounded convex polyhedron. Then
\begin{itemize}
\item[(i)] If $\pol$ has facets with non-parallel unbounded edges, then it can be placed in second trimming position with respect to any of them.
\item[(ii)] If all the unbounded edges of $\pol$ are parallel, then $\pol$ has at least one bounded facet and it can be placed in first trimming position with respect to any of its facets.
\end{itemize}
\end{mprop}

\begin{mprop}[$3$-dimensional polyhedra and second trimming position]\label{2trim}
Let $\pol\subset\R^3$ be a non-degenerate, $3$-dimensional unbounded convex polyhedron. We have:
\begin{itemize}
\item[(i)] If $\pol$ has facets with non-parallel unbounded edges, then it can be placed in second trimming position with respect to any of them.
\item[(ii)] If all the unbounded edges of $\pol$ are parallel, then $\pol$ can be placed in second trimming position with respect to any of its bounded facets.
\end{itemize}
\end{mprop}

\subsection*{Sketch of the proof of the Main Theorem} 
The proof of this result follows from \cite[Thms. 4.1 (ii) \& 6.1 (ii) and Rmks. 4.3 \& 6.3]{fu2} once we 
have guaranteed that: 

\em
\noindent $(*)$ Any non-degenerate $3$-dimensional convex polyhedron $\pol\subset\R^3$ has a \em facet \em $\Ff$ such that $\pol$ can be placed in first and second trimming positions with respect to $\Ff$. 
\em

Of course (*) follows from Propositions \ref{1trim} and \ref{2trim}. Without entering in full detail, the strategy developed in \cite{fu2} (adapted to our $3$-dimensional case) consists of an induction process on the number of facets of the convex polyhedron and it is supported mainly by the following two facts:

(1) A degenerated $3$-dimensional convex polyhedron $\pol\subset\R^3$ can be placed as the product $\p\times\R$ where $\p\subset\R^2$ is a $2$-dimensional convex polygon. In \cite{u2} it is proved that $\R^2\setminus\p$ and $\R^2\setminus\Int\p$ are polynomial images of $\R^2$, say by polynomial maps $f_0,g_0:\R^2\to\R^2$ respectively. Then $\R^3\setminus\pol$ and $\R^3\setminus\Int\pol$ are the images of the polynomial maps $(f_0,\id_\R):\R^3\to\R^3$ and $(g_0,\id_\R):\R^3\to\R^3$ respectively.

(2) Let $\pol:=\{h_1\ge0,\ldots,h_m\ge0\}$ be a non-degenerate $3$-dimensional unbounded convex polyhedron with $m$ facets in (first and second) trimming position with respect to the facet $\Ff$ where each $h_i$ is a linear equation and $\Ff$ is contained in the hyperplane $\{h_1=0\}$. Let $\pol_{\times}:=\{h_2\ge0,\ldots,h_m\ge0\}$ be the $3$-dimensional unbounded convex polyhedron with $m-1$ facets obtained `after eliminating from $\pol$ the facet $\Ff$'. Then there exist by \cite[Lemmas 3.8 \& 5.8]{fu2} polynomial maps $f_1,g_1:\R^3\to\R^3$ such that $f_1(\R^3\setminus\pol_{\times})=\R^3\setminus\pol$ and $g_1(\R^3\setminus\Int\pol_{\times})=\R^3\setminus\Int\pol$ (see \cite[Thms. 4.1 (ii) \& 6.1 (ii) and Rmks 4.3 \& 6.3]{fu2} for the concrete details). Now, by induction hypothesis (and taking into account (1) if $\pol_{\times}$ is degenerated) there exist polynomial maps $f_2,g_2:\R^3\to\R^3$ such that $f_2(\R^3)=\R^3\setminus\pol_{\times}$ and $g_1(\R^3)=\R^3\setminus\Int\pol_{\times}$, so $f:=f_2\circ f_1$ and $g:=g_2\circ g_1$ satisfy the required conditions.
\qed 

\vspace{1mm}
The obstruction to extend the Main Theorem to dimensions $\geq4$ by using our techniques relies on the fact that the following property is exclusive of convex 2-dimensional polyhedra.

\begin{property}\em\label{2sesp}
For any convex polygon $\pol\subset\R^2$ there exists a vectorial line $\vec{\ell}$ and a hyperplane $H$ such that the projection $\pi:\R^2\to H$ with direction $\vec{\ell}$ satisfies the identity $\pi(\pol)=\pi(\pol\cap H)$.
\end{property}

In the Appendix \ref{s4} we exhibit the existence of $n$-dimensional unbounded convex polyhedra that can be placed neither in first trimming position with respect to any of its facets nor in second trimming position with respect to any unbounded facet. 

\subsection*{Structure of the article} The article is organized as follows. All basic notions and (standard) notation appears in Section \ref{s2}. The reading can be started directly in Section \ref{s3} and referred to the preliminaries when needed. In Section \ref{s3} we analyze the boundedness of the set ${\mathfrak A}_\pol$ and provide tools to approach the proofs of Propositions \ref{1trim} and \ref{2trim}, which are developed in Section \ref{s5}. 

\section{Preliminaries on convex polyhedra}\label{s2}

We begin by introducing some preliminary terminology and notations concerning convex polyhedra. For a detailed study of the main properties of convex sets we refer the reader to \cite{ber1,r,z}. An affine hyperplane of $\R^n$ will be denoted as $H:=\{x\in\R^n:\ h(x)=0\}\equiv\{h=0\}$ for a linear equation $h$. It determines two \emph{closed half-spaces}
$$
H^+:=\{x\in\R^n:\ h(x)\geq 0\}\equiv\{h\geq0\}\quad\text{and}\quad H^-:=\{x\in\R^n:\ h(x)\leq 0\}\equiv\{h\leq0\}.
$$
We use an overlying arrow $\vec{\cdot}$ when referring to vectorial staff.

\subsection{Generalities on convex polyhedra}\label{pldrfc} 
A subset $\pol\subset\R^n$ is a \emph{convex polyhedron} if it can be described as the finite intersection $\pol:=\bigcap_{i=1}^rH_i^+$ of closed half-spaces $H_i^+$; we allow this family of half-spaces to by empty to describe $\pol=\R^n$ as a convex polyhedron. The dimension $\dim(\pol)$ of $\pol$ is its dimension as a topological manifold with boundary. By \cite[12.1.5]{ber1} there exists a unique minimal family ${\mathfrak H}:=\{H_1,\ldots,H_m\}$ of affine hyperplanes of $\R^n$, which is empty just in case $\pol=\R^n$, such that $\pol=\bigcap_{i=1}^mH_i^+$; we refer to this family as the \em minimal presentation \em of $\pol$. We assume that we choose the linear equation $h_i$ of each $H_i$ such that $\pol\subset H_i^+$. 

\subsubsection{\!\!}\!\! The \em facets \em or $(n-1)$-\em faces \em of $\pol$ are the intersections $\Ff_i:=H_i\cap\pol$ for $1\leq i\leq m$; only the convex polyhedron $\R^n$ has no facets. Each facet $\Ff_i:=H_i^-\cap\bigcap_{j=1}^mH_j^+$ is a convex polyhedron contained in $H_i$. The convex polyhedron $\pol\subset\R^n$ is a topological manifold with boundary, whose interior is $\Int\pol=\bigcap_{i=1}^m(H_i^+\setminus H_i)$ and its boundary is $\partial\pol=\bigcup_{i=1}^m\Ff_i$. For $0\leq j\leq n-2$ we define inductively the $j$-\em faces \em of $\pol$ as the facets of the $(j+1)$-faces of $\pol$, which are again convex polyhedra. The $0$-faces are the \em vertices \em of $\pol$ and the $1$-faces are the \em edges \em of $\pol$; obviously if $\pol$ has a vertex, then $m\geq n$. A convex polyhedron of $\R^n$ is \em non-degenerate \em if it has at least one vertex; otherwise, we say that the convex polyhedron is \em degenerate\em. 

\subsubsection{\!\!}\!\! A {\em supporting hyperplane} of a convex polyhedron $\pol\subset\R^n$ is a hyperplane $H$ of $\R^n$ that intersects $\pol$ and satisfies $\pol\subset H^+$ or $\pol\subset H^-$. This is equivalent to have $\varnothing\neq\pol\cap H\subset\partial\pol$. The intersection of $\pol$ with a supporting hyperplane $H$ is a face of $\pol$ and conversely each face of $\pol$ is the intersection of $\pol$ with some supporting hyperplane. In particular, the vertices of a convex polyhedron $\pol\subset\R^n$ are those points $p\in\pol$ for which there exists a (supporting) hyperplane $H\subset\R^n$ such that $\pol\cap H=\{p\}$. 

\subsubsection{\!\!}\!\!\label{113} The \em cone of vertex $p$ and base the bounded convex polyhedron $\p\subset\R^n$ \em is the set 
$$
\Cc:=\{\lambda p+(1-\lambda)q: \ q\in\p,\ 0\leq\lambda\leq 1\}.
$$ 
Given $\vec{v}_1,\ldots,\vec{v}_r\in\R^n$ we define the \em cone generated by the vectors $\{\vec{v}_1,\ldots,\vec{v}_r\}$ \em as the set $\vec{\Cc}:=\{\sum_{i=1}^r\lambda_i\vec{v}_i:\ \lambda_i\geq0\}$ and denote $\Cc_p:=p+\vec{\Cc}$ for $p\in\R^n$. Given two points $p,q\in\R^n$ we denote by $\ol{pq}:=\{\lambda p+(1-\lambda)q:\ 0\leq\lambda\leq1\}$ the \em segment connecting $p$ and $q$ \em and given a vector $\vec{v}\in\R^n$, we denote by $p\vec{v}:=\{p+\lambda\vec{v}: \lambda\geq0\}$ the \em half-line of extreme $p$ and direction $\vec{v}$\em. 

\subsection{Recession cone of a convex polyhedron}

We associate to each convex polyhedron $\pol\subset\R^n$ its \em recession cone \em \cite[Ch.1]{z}. Fix a point $p\in\pol$ and denote $\conv{\pol}{}:=\{\vec{v}\in\R^n:\,p\vec{v}\subset\pol\}$\em. Then $\conv{\pol}{}$ is a convex cone and it does not depend on the choice of $p$\em. The set $\conv{\pol}{}$ is called the \em recession cone \em of $\pol$. If $\pol:=\bigcap_{i=1}^rH_i^+$, then $\conv{\pol}{}:=\bigcap_{i=1}^r\conv{H_i^+}{}=\bigcap_{i=1}^r\vec{H_i}^+$. Clearly, $\conv{\pol}{}=\{{\bf0}\}$ if and only if $\pol$ is bounded. In addition, if $\p\subset\R^n$ is a non-degenerate convex polyhedron and $k\geq1$, then $\conv{\R^k\times\p}{}=\R^k\times\conv{\p}{}$. Recall that each degenerate convex polyhedra can be written as the product of a non-degenerate convex polyhedron times an Euclidean space.

\subsubsection{\!\!}\!\! \label{conoedges}
An important property of a bounded convex polyhedron is that it coincides with the \em convex hull \em of the set of its vertices. A general non-degenerate convex polyhedron $\pol\subset\R^n$ can be described as follows \cite[Ch.1]{z}. \em Let ${\mathfrak V}:=\{p_1,\ldots,p_r\}$ be the set of vertices of $\pol$ and let ${\mathfrak A}:=\{\Aa_1,\ldots,\Aa_s\}$ be the set of unbounded edges of $\pol$. Write $\Aa_j:=q_j\vec{v}_j$ for $j=1,\ldots,s$. Then 
\begin{itemize}
\item[(i)] $\conv{\pol}{}=\{\sum_{j=1}^s\lambda_j\vec{v}_j:\ \lambda_1,\ldots,\lambda_s\geq0\}$ or $\{{\bf0}\}$ if ${\mathfrak A}=\varnothing$.
\item[(ii)] $\pol=\pol_0+\conv{\pol}{}$ where $\pol_0$ is the bounded convex polyhedron of vertices $p_1,\ldots,p_r$.
\end{itemize}\em

\subsection{Facing upwards positions for convex polyhedra}

The proof of the Main Theorem has been reduced to show that non-degenerate unbounded convex polyhedra of $\R^3$ can be placed in a specific form. As an initial step we introduce the concept of \em facing upwards positions\em.

\begin{define}\label{def:upwards}
An unbounded non-degenerate convex polyhedron $\pol\subset\R^n$ is in \em facing upwards position with respect to the hyperplane $\vec{\Pi}$ \em of $\R^n$ (shortly, {\rm FU}-position w.r.t. $\vec{\Pi}$) if there exist a hyperplane $\Pi$ parallel to $\vec{\Pi}$ that intersects all the unbounded edges of $\pol$ and such that all the vertices of $\pol$ belong to the open half-space $\Int\Pi^-$. The hyperplane $\Pi$ is called a \em sawing hyperplane \em for $\pol$. Any hyperplane $\Pi'\subset\Pi^+$ (parallel to $\vec{\Pi}$) is also a sawing hyperplane for $\pol$.

Let $\vec{h}$ be a linear equation of $\vec{\Pi}$. We say that the FU-position of $\pol$ with respect to $\vec{\Pi}$ is \em optimal \em if the minimum of $\vec{h}|_{\pol}$ is attained exactly in one point (which must be a vertex of $\pol$). 
\end{define}

\subsubsection{Connection with bounded convex polyhedra.}\label{lem:bounded} 
Each unbounded non-degenerate convex polyhedron $\pol\subset\R^n$ can be placed in optimal FU-position w.r.t. the hyperplane $\vec{\Pi}:=\{x_n=0\}$ in such a way it does not intersect the hyperplane $\Pi_0:=\{x_n=0\}$. Under this hypothesis there exists a natural bridge between non-degenerate unbounded convex polyhedra and bounded ones. Denote the hyperplane at infinity of the projective space $\R\PP^n$ with $\mathsf{H}_\infty(\R)$. Write $\widehat{\pol}:=\cl_{\R\PP^n}(\pol)=\pol\sqcup\pol_\infty$ where $\pol_\infty:=\cl_{\R\PP^n}(\pol)\cap\mathsf{H}_\infty(\R)$ and consider the involution
$$
\phi:\R\PP^n\to\R\PP^n,\ (x_0:x_1:\ldots:x_{n-1}:x_n)\mapsto(x_n:x_1:\ldots:x_{n-1}:x_0),
$$
induced by the birational map
$$
f:=\phi|_{\R^n}:\R^n\dashrightarrow\R^n,\ (x_1,\ldots,x_n)\mapsto(y_1,\ldots,y_n):=\Big(\frac{x_1}{x_n},\ldots,\frac{x_{n-1}}{x_n},\frac{1}{x_n}\Big).
$$
Then \em $\pol':=\cl(f(\pol))=\phi(\widehat{\pol})\subset\R^n\equiv\{y_0\neq0\}$ is a bounded convex polyhedron and one of its faces is $\Ee':=\phi(\pol_\infty)$. Moreover $\phi(p+\conv{\pol}{})\cup\phi(\pol_\infty)$ is the closed cone $\Cc_{\phi(p)}$ of base $\phi(\pol_\infty)$ and vertex $\phi(p)$ for each $p\in\pol$\em.
\begin{proof}
We just check that $\phi(\widehat{\pol})\subset\R^n$. Otherwise, as $\pol\cap\Pi_0=\varnothing$, there exists $z:={(0:\vec{v})}\in\pol_\infty\cap\cl_{\R\PP^n}(\Pi_0)$ where $\vec{v}:=(v_1,\ldots,v_{n-1},0)\in\vec{\Pi}$. Since $\pol$ is in optimal FU-position w.r.t. $\{x_n=0\}$, there exists a vertex $p:=(p_1,\ldots,p_n)$ of $\pol$ such that $p_n<x_n$ for all $x:=(x_1,\ldots,x_n)\in\pol\setminus\{p\}$. Let $\{z_k\}_k\subset\pol\setminus\{p\}$ be a sequence that converges to $z$. For each $k\geq1$ choose a unitary vector $\vec{u}_k$ and $\lambda_k>0$ such that $z_k=p+\lambda_k\vec{u}_k$. By the compactness of the unitary sphere $\sph^{n-1}$ we assume that the sequence $\{\vec{u}_k\}_k$ converges to a unitary vector $\vec{u}$. Note that $\lim_{k\to\infty}\lambda_k=+\infty$ because $z\in\pol_\infty$. Observe that $z_k\equiv(1:z_k)=(\frac{1}{\lambda_k}:\frac{p}{\lambda_k}+\vec{u}_k)$ tends to $(0:\vec{u})$ when $k\to+\infty$, so assume $\vec{u}=\vec{v}$ and $z:=(0:\vec{u})$. We see now that the half-line $p\vec{u}\subset\pol$. Indeed, let $x\in p\vec{u}$ and write $x=p+\|x-p\|\vec{u}$. Let $k_0\geq1$ be such that for each $k\geq k_0$ it holds $\lambda_k>\|x-p\|$; hence, 
$$
y_k=p+\|x-p\|\vec{u}_k\in\ol{pz_k}=\{p+t\vec{u}_k:\ 0\leq t\leq\lambda_k\}\subset\pol.
$$
As $\pol$ is closed, $x:=p+\|x-p\|\vec{u}=\lim_{k\to\infty}y_k\in\pol$. Thus, 
$$
\{(p_1+tv_1,\ldots,p_{n-1}+tv_{n-1},p_n):\ t\geq0\}=p\vec{u}\subset\pol,
$$ 
against the fact that $p_n<x_n$ for all $x:=(x_1,\ldots,x_n)\in\pol\setminus\{p\}$.
\end{proof}

\subsubsection{Supporting hyperplanes} There is an easy procedure to find supporting hyperplanes for a non-degenerate unbounded convex polyhedron $\pol\subset\R^n$ in FU-position w.r.t. a hyperplane $\vec{\Pi}$. 

\begin{lem}\label{ppdrhs}
Let $\Pi$ be a sawing hyperplane for $\pol$ and let $W$ be a supporting hyperplane in $\Pi$ of the convex polyhedron $\p:=\pol\cap\Pi$. Let $p\in W\cap\p$ be a vertex of $\p$ and let $\Aa$ be the unbounded edge of $\pol$ such that $\{p\}=\Aa\cap\Pi$. Then the affine subspace $H$ of $\R^n$ generated by $W$ and $\Aa$ is a supporting hyperplane of $\pol$. 
\end{lem}

\subsubsection{Recession cone with maximal dimension} We finish with a technical result concerning the recession cone $\conv{\pol}{}$ of an $n$-dimensional non-degenerate unbounded convex polyhedron $\pol\subset\R^n$ in \em FU-\em position w.r.t. $\vec{\Pi}$ when $\dim(\conv{\pol}{})=n$.

\begin{lem}\label{lem:cono4} 
Assume $\dim(\conv{\pol}{})=n$, let $\vec{v}\in\Int\conv{\pol}{}$ and consider a finite set ${\mathfrak G}\subset\R^n$. Then there exists a hyperplane $\Pi$ parallel to $\vec{\Pi}$ such that $p\vec{v}\cap(\pol\cap\Pi)$ is a singleton for each $p\in{\mathfrak G}$.
\end{lem}
\begin{proof} 
We use freely the straightforward fact: 

\vspace{1mm}\setcounter{substep}{0}
\begin{substeps}{lem:cono4}\label{stf}\em
For each $p\in\R^n$ the intersection $p\vec{v}\cap\pol$ is a half-line $p_1\vec{v}\subset\pol$.
\end{substeps}

\vspace{1mm}
Write $\vec{\Pi}:=\{\vec{h}=0\}$ and let $\Pi_0:=\{h_0:=a_0+\vec{h}=0\}$ be a sawing hyperplane for $\pol$. Write $\Aa_i:=q_i\vec{v}_i$ where $q_i$ is a vertex of $\pol$, the vector $\vec{v}_i\in\conv{\pol}{}$ and $b_i:=q_i+\vec{v}_i\in\Pi$ for $i=1,\ldots,s$. As $0=h_0(b_i)=h_0(q_i)+\vec{h}(\vec{v}_i)$ and $h_0(q_i)<0$, we deduce that $\vec{h}(\vec{v}_i)=-h_0(q_i)>0$. By \ref{conoedges} there exist $\lambda_i\geq0$ not all zero such that $\vec{v}=\sum_{i=1}^s\lambda_i\vec{v}_i$; hence,
$$
\vec{h}(\vec{v})=\vec{h}\Big(\sum_{i=1}^s\lambda_i\vec{v}_i\Big)=\sum_{i=1}^s\lambda_i\vec{h}(\vec{v}_i)>0.
$$
For each $p\in{\mathfrak G}$ let $p'\in\pol$ be such that $p\vec{v}\cap\pol=p'\vec{v}$. We choose $\Pi=\{h=0\}$ parallel to $\vec{\Pi}$ such that $p'\in\Int{\Pi'}^-$ for each $p\in{\mathfrak G}$. We claim that: \em $p\vec{v}\cap(\pol\cap\Pi')\neq\varnothing$ for each $p\in{\mathfrak G}$\em. Indeed, fix $p\in{\mathfrak G}$. As $\vec{h}(\vec{v})>0$ and $h(p')<0$, there exists $t>0$ such that $h(p'+t\vec{v})=0$, so $p'+t\vec{v}\in p'\vec{v}\cap\Pi=p\vec{v}\cap(\pol\cap\Pi)$. Obviously, as $p'\vec{v}$ is a half-line and $p'\not\in\Pi$, the intersection $p\vec{v}\cap(\pol\cap\Pi)$ is a singleton, as required.
\end{proof}

\section{Characterization of the boundedness of the set ${\mathfrak A}_\pol$}\label{s3}

In this section we characterize when the set ${\mathfrak A}_\pol$ of a non-degenerate unbounded convex polyhedron $\pol$ is bounded. Recall that
$$
\pi_n:\R^n\to\R^{n-1}\times\{0\},\ x:=(x_1,\ldots,x_n)\mapsto(x',0):=(x_1,\ldots,x_{n-1},0)
$$
and we denote by $\vec{\ell}_n$ the vectorial line generated by the vector $\vec{e}_n:=(0,\ldots,0,1)$. 

\begin{lem}[Boundedness of ${\mathfrak A}_\pol$]\label{proysec}
Let $\pol\subset\R^n$ be a non-degenerate unbounded convex polyhedron. Then the set ${\mathfrak A}_\pol:=\{a\in\R^{n-1}:\ \Ii_a:=\pi_n^{-1}(a,0)\cap\pol\neq\varnothing,\ (a,0)\not\in\Ii_a\}$ is bounded if and only if the following conditions hold:
\begin{itemize}
\item[(i)] $\vec{\pi}_n(\conv{\pol}{})=\conv{\pol}{}\cap\{x_n=0\}$.
\item[(ii)] There exists a hyperplane $\Pi\subset\R^n$ parallel to $\vec{\ell}_n$ such that: it meets all the unbounded edges of $\pol$ that are non-parallel to $\vec{\ell}_n$, all the vertices of $\pol$ and all the unbounded edges of $\pol$ parallel to $\vec{\ell}_n$ are contained in $\Int\Pi^-$ and $\pi_n(\pol\cap\Pi)=(\pol\cap\Pi)\cap\{x_n=0\}$.
\end{itemize}
\end{lem}

\subsection{Hyperplane sections of a convex polyhedron}
In order to prove Lemma \ref{proysec} we need to understand the generic sections parallel to a hyperplane $\vec{\Pi}:=\{\vec{h}=0\}$ of an $n$-dimen\-sional non-degenerate unbounded convex polyhedron $\pol\subset\R^n$. 

\begin{lem}\label{sectioning0} 
Let $\Aa_1,\ldots,\Aa_s$ be the unbounded edges of $\pol$ and assume that the first $k$ are non-parallel to $\vec{\Pi}$ and the remaining $s-k$ are parallel to $\vec{\Pi}$. For each $i=1,\ldots,s$ write $\Aa_i:=q_i\vec{v}_i$ where $q_i$ is a vertex of $\pol$ and $\vec{v}_i\in\conv{\pol}{}$. Let $\Pi_0$ be a hyperplane parallel to $\vec{\Pi}$ that meets $\Aa_1,\ldots,\Aa_k$ and such that all the vertices of $\pol$ and the unbounded edges $\Aa_{k+1},\ldots,\Aa_s$ are contained in $\Int\Pi_0^-$. For each hyperplane $\Pi\subset\Pi_0^+$, we have: 
\begin{itemize}\em
\item[(i)] $\p:=\pol\cap\Pi$ is the convex polyhedron whose vertices are the intersections of $\Pi$ with $\Aa_1,\ldots,\Aa_k$ and its recession cone is $\conv{\p}{}=\conv{\pol\cap{\Pi}^-}{}=\{\sum_{j=k+1}^s\lambda_j\vec{v}_j:\ \lambda_j\geq0\}$.
\item[(ii)] $\pol\cap\Pi^+=\p+\conv{\pol}{}=\p+\{\sum_{j=1}^k\lambda_j\vec{v}_j:\ \lambda_j\geq0\}$.
\end{itemize}
\end{lem}
\begin{remark}\label{sectioning}
If $k=s$, then $\pol$ is in FU-position w.r.t. $\vec{\Pi}$ and it holds:
\begin{itemize}
\item[(i)] $\p:=\pol\cap\Pi\neq\varnothing$ is bounded and its vertices are the intersections of $\Pi$ with $\Aa_1,\ldots,\Aa_s$.
\item[(ii)] $\pol\cap\Pi^-$ is bounded and $\pol\cap\Pi^+=\p+\conv{\pol}{}$.
\end{itemize}
\end{remark}
\begin{proof}[Proof of Lemma \em\ref{sectioning}]
(i) Let $p$ be a vertex of $\p$. As $p\in\partial\pol$, we choose a face $\Ee$ of $\pol$ of the smallest dimension between those containing $p$; clearly, $p\in\Int\Ee$. As the vertices of $\pol$ are contained in $\Int\Pi^-$ and $\Pi\subset\Pi^+$, the dimension of $\Ee$ is $\geq1$. Let $W$ be the affine subspace generated by $\Ee$ and notice that $W\cap\Pi=\{p\}$ (otherwise we would have a face of $\p$ crossing the vertex $p$). As $\Pi$ is a hyperplane, $\dim W=1$. Thus, $\Ee:=\Aa_j$ is an (unbounded) edge of $\pol$ non-parallel to $\vec{\Pi}$. On the other hand, it is clear that the intersections of $\Pi$ with $\Aa_1,\ldots,\Aa_k$ are vertices of $\p$.

Next, we prove that $\conv{\p}{}=\conv{\pol\cap{\Pi}^-}{}=\vec{\Cc}:=\{\sum_{j=k+1}^s\lambda_j\vec{v}_j:\ \lambda_j\geq0\}$.

Indeed, let $h:=a+\vec{h}$ be a linear equation of $\Pi$ such that $\Pi^+=\{h\geq0\}$. It holds 
\begin{equation}\label{signoshv}
\vec{h}(\vec{v}_i)\begin{cases} 
>0&\text{for $i=1,\ldots,k$,}\\
=0&\text{for $i=k+1,\ldots,s$.} 
\end{cases}
\end{equation}
If $\vec{v}\in\conv{\pol\cap{\Pi}^-}{}$, then $\vec{v}\in\conv{\pol}{}$ and $\vec{h}(\vec{v})\leq0$. As $\conv{\pol}{}=\{\sum_{i=1}^s\lambda_i\vec{v}_i:\ \lambda_i\geq0\}$, there exist $\zeta_i\geq0$ such that $\vec{v}=\sum_{j=1}^s\zeta_j\vec{v}_j$. By \eqref{signoshv}
$$
0\geq\vec{h}(\vec{v})=\vec{h}\Big(\sum_{i=1}^s\zeta_i\vec{v}_i\Big)=\sum_{i=1}^s\zeta_i\vec{h}(\vec{v}_i)=\sum_{i=1}^k\zeta_i\vec{h}(\vec{v}_i)\geq0;
$$
hence, $\zeta_j=0$ for $j=1,\ldots,k$, so $\vec{v}\in\vec{\Cc}$. Consequently,
$$
\conv{\pol\cap{\Pi}^-}{}\subset\vec{\Cc}\subset\conv{\pol}{}\cap\vec{\Pi}=\conv{\pol\cap\Pi}{}=\conv{\p}{}\subset\conv{\pol\cap{\Pi}^-}{}.
$$

(ii) The vertices of $\pol\cap{\Pi}^+$ are the intersections of $\Pi$ with the edges $\Aa_1,\ldots,\Aa_k$ of $\pol$. Thus, the convex hull of the set consisting of those vertices is contained in $\p$, so $\pol\cap{\Pi}^+=\p+\conv{\pol\cap{\Pi}^+}{}$. As $\conv{\pol}{}=\conv{\pol\cap{\Pi}^+}{}$ and $\conv{\p}{}=\{\sum_{j=k+1}^s\lambda_j\vec{v}_j:\ \lambda_j\geq0\}$, we deduce
$$
\pol\cap{\Pi}^+=\p+\conv{\pol}{}=\p+\conv{\p}{}+\Big\{\sum_{j=1}^k\lambda_j\vec{v}_j:\ \lambda_j\geq0\Big\}=\p+\Big\{\sum_{j=1}^k\lambda_j\vec{v}_j:\ \lambda_j\geq0\Big\},
$$
as required.
\end{proof}

\subsection{Proof of Lemma \ref{proysec}} Denote the vertices of $\pol$ with $p_1,\ldots,p_r$ and the unbounded edges of $\pol$ with $\Aa_1,\ldots,\Aa_s$ ordered in such a way that the first $k$ are non-parallel to $\vec{\ell}_n$ and the remaining $s-k$ are parallel to $\vec{\ell}_n$. For each $i=1,\ldots,s$ we write $\Aa_i:=q_i\vec{v}_i$ where $q_i$ is a vertex of $\pol$ and $\vec{v}_i\in\conv{\pol}{}$ for $i=1,\ldots,s$. 

Assume that ${\mathfrak A}_\pol$ is bounded. We prove first (i) $\vec{\pi}_n(\conv{\pol}{})=\conv{\pol}{}\cap\{x_n=0\}$. Write $H:=\{x_n=0\}$ and let us see that if $\vec{v}:=(v_1,\ldots,v_{n})\in\conv{\pol}{}$, then 
$$
\vec{w}:=\vec{\pi}_n(\vec{v})=(v_1,\ldots,v_{n-1},0)\in\conv{\pol}{}\cap\vec{H}.
$$ 

Indeed, if $\vec{v}\in\vec{\ell}_n$, then $\vec{\pi}_n(\vec{v})={\bf0}\in\conv{\pol}{}\cap\vec{H}$, so assume $\vec{v}\not\in\vec{\ell}_n$. Let $p\in\pol$; as $\vec{v}\in\conv{\pol}{}$, the half-line $p\vec{v}\subset\pol$. In addition $\vec{v}\not\in\vec{\ell}_n$, so $\vec{w}\in\vec{H}\setminus\{{\bf0}\}$. As the set ${\mathfrak A}_\pol$ is bounded and $\pi_n(p\vec{v})=\pi_n(p)\vec{w}$ is a half-line, there exists a point $q\in p\vec{v}$ such that the half-line $\pi_n(q\vec{v})=\pi_n(q)\vec{w}$ does not meet ${\mathfrak A}_\pol\times\{0\}\subset H$. For each $t>0$ the point $q+t\vec{v}\in\pol\cap\pi_n^{-1}(\pi_n(q)+t\vec{w})$ while $\pi_n(q)+t\vec{w}\not\in{\mathfrak A}_\pol\times\{0\}$; hence, the half-line $\pi_n(q)\vec{w}\subset\pol\cap H$ and $\vec{w}\in\conv{\pol}{}\cap\vec{H}$. Thus, 
$$
\vec{\pi}_n(\conv{\pol}{})\subset\conv{\pol}{}\cap\vec{H}
$$
and consequently $\vec{\pi}_n(\conv{\pol}{})=\conv{\pol}{}\cap\{x_n=0\}$ because the other inclusion is trivial. 

To show (ii) we distinguish two cases:

\noindent{\sc Case 1}. If $k=0$, the edges $\Aa_1,\ldots,\Aa_s$ are parallel to $\vec{\ell}_n$. Then $\pol$ is in FU-position w.r.t. $\{x_n=0\}$ and by Remark \ref{sectioning} the projection $\pi_n(\pol)$ is bounded. We choose a hyperplane $\Pi$ parallel to the line $\vec{\ell}_n$ such that $\pol\subset\Int\Pi^-$. Clearly, $\Pi$ enjoys the required conditions.

\noindent{\sc Case 2}. If $k>0$, we know that $\conv{\pol}{}=\{\sum_{i=1}^s\lambda_i\vec{v}_i:\ \lambda_i\geq0\}$ and, if $\vec{w}_i:=\vec{\pi}_n(\vec{v}_i)$, we deduce 
$$
\vec{\pi}_n(\conv{\pol}{})=\Big\{\sum_{i=1}^s\lambda_i\vec{w}_i:\ \lambda_i\geq0\Big\}=\Big\{\sum_{i=1}^k\lambda_i\vec{w}_i:\ \lambda_i\geq0\Big\}\neq\{{\bf0}\}
$$
is an unbounded cone whose unique vertex is ${\bf0}$ (otherwise $\pol$ would be degenerate). Observe that $\vec{w}_i\neq0$ for each $i=1,\ldots,k$, because the edges $\Aa_1,\ldots,\Aa_k$ are non-parallel to $\vec{\ell}_n$.

Let $W$ be a supporting hyperplane of $\vec{\pi}_n(\conv{\pol}{})$ in $\{x_n=0\}$ such that $\vec{\pi}_n(\conv{\pol}{})\cap W=\{{\bf0}\}$ and consider the hyperplane $\vec{\Pi}:=\vec{W}+\vec{\ell}_n=\{\vec{h}=0\}$. We may assume that $\vec{h}(\vec{w}_i)>0$ for $i=1,\ldots,k$. As $\vec{v}_i-\vec{w}_i\in\vec{\ell}_n+t_i$, we have $\vec{h}(\vec{v_i})=\vec{h}(\vec{w}_i)>0$.

Let $\Pi:=\{h=0\}$ be a hyperplane parallel to $\vec{\Pi}$ such that all the vertices $p_i$ of $\pol$ and the set ${\mathfrak A}_\pol\times\{0\}$ are contained in $\Int\Pi^-$. As $\vec{h}(\vec{v_i})=\vec{h}_0(\vec{v_i})>0$, it holds that $\Pi$ meets the edges $\Aa_i=p_i\vec{v}_i$ for $i=1,\ldots,k$. Since $\vec{\ell}_n\subset\vec{\Pi}$,
$$
\pi_n^{-1}({\mathfrak A}_\pol\times\{0\})\subset\pi_n^{-1}(\pi_n(\{h<0\}))=\{h<0\}=\Int\Pi^-;
$$
hence, $\pi_n(x)\in\pol\cap\Pi\cap\{x_n=0\}$ for each $x\in\pol\cap\Pi$. Thus, $\pi_n(\pol\cap\Pi)=(\pol\cap\Pi)\cap\{x_n=0\}$.

\vspace{1mm}
Assume now that (i) and (ii) hold. We distinguish two cases to prove that ${\mathfrak A}_\pol$ is bounded:

\noindent{\sc Case 1.} $\pol\cap\Pi=\varnothing$, o equivalently, all the unbounded edges of $\pol$ are parallel to the line $\vec{\ell}_n$. Then $\pol$ is in FU-position w.r.t. $\{x_n=0\}$ and by Remark \ref{sectioning} the projection $\pi_n(\pol)$ is a bounded convex polyhedron. Thus, ${\mathfrak A}_\pol\times\{0\}\subset\pi_n(\pol)$ is a bounded set.

\noindent{\sc Case 2.} $\pol\cap\Pi\neq\varnothing$. By hypothesis (ii) and Lemma \ref{sectioning0}
\begin{equation}\label{poc}
\pol\cap{\Pi}^+=(\pol\cap\Pi)+\conv{\pol}{}\ \text{ and }\ \conv{\pol\cap{\Pi}^-}{}=\Big\{\sum_{j=k+1}^s\lambda_j\vec{v}_j:\ \lambda_j\geq0\Big\}.
\end{equation}
The vertices of $\pol\cap{\Pi}^-$ are the vertices of $\pol$ and the intersections $\{y_j\}=\Aa_j\cap\Pi$ for $j=1,\ldots,k$. By \ref{conoedges} $\pol\cap\Pi^-=\pol_0+\conv{\pol\cap{\Pi}^-}{}$ where $\pol_0$ is the convex hull of $\{p_1,\ldots,p_r,y_1,\ldots,y_k\}$. As $\vec{\pi}_n(\vec{v}_j)={\bf0}$ for $j=k+1,\ldots,s$,
$$
\pi_n(\pol\cap{\Pi}^-)=\pi_n(\pol_0+\conv{\pol\cap{\Pi}^-}{})=\pi_n(\pol_0)+\vec{\pi}_n(\conv{\pol\cap{\Pi}^-}{})=\pi_n(\pol_0),
$$
which is a bounded set. Thus, if we prove that $\pi_n(\pol\cap{\Pi}^+)=\pol\cap{\Pi}^+\cap\{x_n=0\}$, the set ${\mathfrak A}_\pol\times\{0\}\subset\pi_n(\pol_0)$ is bounded. Indeed, by hypotheses (i), (ii) and equation \eqref{poc} it holds
$$
\pi_n(\pol\cap{\Pi}^+)=\pi_n(\pol\cap\Pi)+\vec{\pi}_n(\conv{\pol}{})=((\pol\cap\Pi)+\conv{\pol}{})\cap\{x_n=0\}=\pol\cap{\Pi}^+\cap\{x_n=0\},
$$
as required.
\qed

\section{Placing in trimming positions}\label{s5}

\subsection{First trimming position}

The purpose of the first part of this section is to prove Proposition \ref{1trim}. For the sake of clearness, we divide the proof into two parts.

\begin{proof}[Proof of Proposition \em \ref{1trim} \em for unbounded facets]
Assume $\pol$ is in FU-position w.r.t. $\vec{\Pi}=\{x_n=0\}$. The proof runs in several steps (see Figure~\ref{prirec}):

\vspace{1mm}
\noindent{\bf Step 1}.\,\em Choice of a suitable facet $\Ff$\em. Denote $\widehat{\pol}:=\cl_{\R\PP^3}(\pol)$ and $\pol_\infty:=\cl_{\R\PP^3}(\pol)\cap\mathsf{H}_\infty(\R)$. Consider the involutive homography
$$
\phi:\R\PP^3\to\R\PP^3,\ (x_0:x_1:x_2:x_3)\mapsto(x_3:x_1:x_2:x_0)=(y_0:y_1:y_2:y_3).
$$
By \ref{lem:bounded}, the convex polyhedron $\pol':=\phi(\widehat{\pol})\subset\R^3\equiv\{y_0\neq0\}$ is bounded and one of its faces is $\Ee':=\phi(\pol_\infty)$, which needs not to be a facet. Now we distinguish:

\noindent{\sc Case 1.} If $\Ee'$ has dimension $0$, choose any unbounded facet $\Ff$ of $\pol$.

\noindent{\sc Case 2.} If $\Ee'$ has dimension $1$ ($\Ee'$ is an edge of $\phi(\widehat{\pol})$) or $2$ ($\Ee'$ is a facet of $\phi(\widehat{\pol})$), choose any unbounded facet $\Ff$ with non-parallel unbounded edges. Note that $\dim(\cl(\phi(\Ff))\cap\Ee')=1$. 

Denote by $H$ the plane generated by the facet $\Ff$.

\vspace{1mm}
\noindent{\bf Step 2.} \em Choice of a suitable sawing plane $\Pi$ and an auxiliary vector $\vec{w}$\em.
Let $\Pi$ be sawing plane for $\pol$ (parallel to $\vec{\Pi}$). If $\dim(\Ee')=2$, the dimension of $\pol_\infty$ is $2$, the dimension of $\cl_{\R\PP^3}(\Ff)\cap\pol_\infty$ is $1$ and by \ref{lem:bounded} the dimension of $\conv{\Ff}{}$ is $2$. Since $\pol$ is in FU-position w.r.t. $\vec{\Pi}$, then $\Ff$ is in FU-position w.r.t. $\vec{\Pi}\cap\vec{H}$ inside $H$. Let ${\mathfrak G}$ be the (finite) set constituted by the intersections of $H$ with the unbounded edges of $\pol$ that are non-parallel to $H$ and fix $\vec{w}\in\Int\conv{\Ff}{}$. By Lemma \ref{lem:cono4} we may assume in addition that the intersection $p\vec{w}\cap(\Ff\cap\Pi)$ is a singleton for each $p\in{\mathfrak G}$.

\vspace{1mm}
\noindent{\bf Step 3.} \em Construction of an auxiliary supporting hyperplane $H_0$ of $\pol$\em.
Denote $\p:=\pol\cap\Pi$ and note that $\Ff\cap\Pi=H\cap\p$ is one of its edges. The map $\rho:\Pi\to\Pi/(\vec{H}\cap\vec{\Pi}),\ x\mapsto x+(\vec{H}\cap\vec{\Pi})$ is continuous (with respect to the quotient topology of $\Pi/(\vec{H}\cap\vec{\Pi})$). As $\p\subset\Pi$ is compact and $\Pi/(\vec{H}\cap\vec{\Pi})$ is homeomorphic to $\R$ (with its usual topology), $\rho(\p)\equiv[a,b]$ is a compact interval (nontrivial because $\p$ has dimension $2$). Clearly, $\rho(H\cap\Pi)$ is one of the extremes of this interval and we assume that $\Lambda_a:=\rho^{-1}(a)=H\cap\Pi$. Note that $\Lambda_b:=\rho^{-1}(b)$ is a supporting line of $\p$ in $\Pi$, so $\Lambda_b\cap\p$ is either an edge or a vertex of $\p$. We pick a vertex $p_0\in\Lambda_b\cap\p$, which is the intersection of $\Pi$ with an unbounded edge $\Aa$ of $\pol$. Denote the line generated by $\Aa$ with $\ell$.

Let $H_0$ be the supporting plane of $\pol$ generated by $\Lambda_b$ and $\ell$ (Lemma \ref{ppdrhs}). We can assume that $\pol\subset H^+\cap H_0^+$.

\vspace{1mm}
\noindent{\bf Step 4.}
\em Construction of a plane $W$ that contains the line $\ell$ and such that $W\cap\Ff$ is a half-line\em. In order to achieve this we analyze two possible situations:

\noindent{\sc Case 1.} If $\ell$ is parallel to $H$, choose a plane $W$ that contains $\ell$ and meets $\Int\Ff\cap\Pi$ in a point $q$. 

\noindent{\sc Case 2.} If $\ell$ is non-parallel to $H$, then it meets $H$ in a singleton, say $\ell\cap H=\{q_0\}$. Besides, 
$$
\cl_{\R\PP^3}(\ell)\cap\cl_{\R\PP^3}(H)\cap\mathsf{H_\infty}(\R)=\varnothing\quad \text{implies}\quad \cl_{\R\PP^3}(\Aa)\cap\cl_{\R\PP^3}(\Ff)\cap\pol_\infty=\varnothing,
$$ 
so the dimension of $\pol_\infty$ is $2$. Recall now the vector auxiliary $\vec{w}$ fixed in Step 2. By our choice of $\Pi$ (see Step 2) the intersection $q_0\vec{w}\cap(\Ff\cap\Pi)$ is a point $q$. Denote the plane of $\R^3$ that contains the coplanar lines $\ell$ and the line through $q_0$ that is parallel to $\vec{w}$ with $W$.

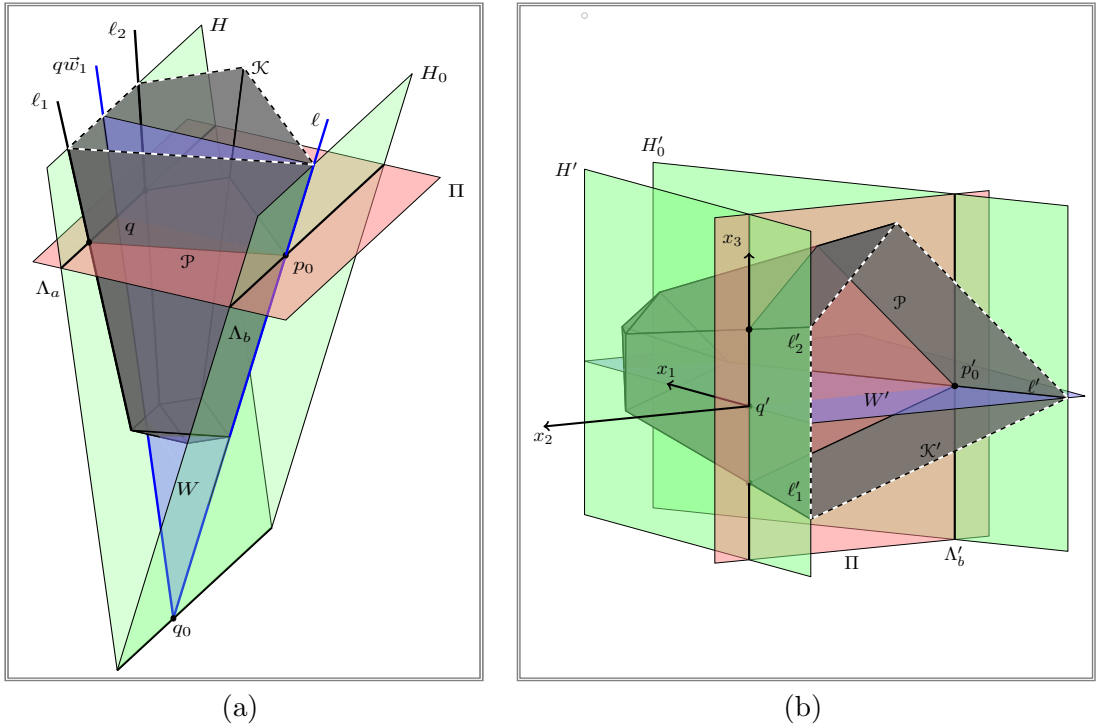
\begin{figure}[ht]\label{fig3}
\begin{center}
\begin{tabular}{cc}
\resizebox{6.4cm}{9cm}{
\begin{tikzpicture}[line join=round,background rectangle/.style=
{double, thick, draw=gray},
show background rectangle, xscale=0.6, yscale=0.8]
\tikzstyle{planelinestyle} = [cull=false,red,dashed]
\tikzstyle{blinestyle} = [cull=false,red,dashed]
\tikzstyle{plinestyle} = [cull=false,red,dashed]
\tikzstyle{verstyle}=[cull=false,fill=blue!20,fill opacity=0.8]
\tikzstyle{horstyle}=[cull=false,fill=blue!20,fill opacity=0.8]
\tikzstyle{boxstyle}=[cull=false,fill=blue!20,fill opacity=0.8]
\tikzstyle{polystyle} = [cull=false,fill=blue!20,fill opacity=0.8]
\draw[very thick](.408,-4.667)--(4.899,-1);
\filldraw[fill=red!50,fill opacity=0.5,draw=none](1.225,7.667)--(0,7.5)--(2.449,9.5)--(3.266,9.333)--cycle;
\draw[draw=black,thin](0,7.5)--(2.449,9.5)--(3.266,9.333);
\filldraw[fill=red!50,fill opacity=0.5,draw=none](1.225,7.667)--(.408,7)--(-.408,7.167)--(0,7.5)--cycle;
\draw[draw=black,thin](-.408,7.167)--(0,7.5);
\filldraw[fill=red!50,fill opacity=0.5,draw=none](.408,7)--(-.408,6.333)--(-1.361,6.389)--(-.408,7.167)--cycle;
\draw[draw=black,thin](-1.361,6.389)--(-.408,7.167);
\filldraw[fill=red!50,fill opacity=0.5,draw=none](-.408,6.333)--(-1.225,5.667)--(-2.041,5.833)--(-1.361,6.389)--cycle;
\draw[draw=black,thin](-1.225,5.667)--(-2.041,5.833)--(-1.361,6.389);
\filldraw[draw=black,thin,fill=green!50,fill opacity=0.3](-1.633,8.25)--(2.858,11.917)--(4.899,-1)--(.408,-4.667)--cycle;
\filldraw[fill=red!50,fill opacity=0.5,draw=none](3.674,8)--(2.858,9)--(3.266,9.333)--(7.348,8.5)--cycle;
\draw[draw=black,thin](3.266,9.333)--(7.348,8.5);
\draw[very thick](1.225,7.667)--(3.266,9.333);
\filldraw[fill=red!50,fill opacity=0.5,draw=none](3.674,8)--(1.225,7.667)--(2.858,9)--cycle;
\filldraw[draw=black,thick,fill=black!60,fill opacity=0.8](1.021,10.417)--(4.082,10.833)--(2.858,2.333)--(1.633,2.167)--cycle;
\filldraw[fill=red!50,fill opacity=0.5,draw=none](5.307,6)--(3.674,8)--(7.348,8.5)--(8.165,8.333)--cycle;
\draw[draw=black,thin](7.348,8.5)--(8.165,8.333);
\filldraw[fill=black!60,fill opacity=0.8,draw=none](3.674,8)--(5.307,6)--(3.674,1.333)--(2.858,2.333)--cycle;
\draw[draw=black,thick](5.307,6)--(3.674,1.333)--(2.858,2.333)--(3.674,8);
\filldraw[draw=black,thick,fill=black!60,fill opacity=0.8](2.449,1.167)--(1.633,2.167)--(2.858,2.333)--cycle;
\filldraw[draw=black,thick,fill=black!60,fill opacity=0.8](2.449,1.167)--(2.858,2.333)--(3.674,1.333)--cycle;
\filldraw[fill=blue!50,fill opacity=0.5,draw=none](1.633,1.333)--(1.304,1.333)--(2.041,-3.333)--(3.062,-.417)--cycle;
\draw[draw=black,thin](1.304,1.333)--(2.041,-3.333)--(3.062,-.417);
\draw[very thick,draw=blue](1.304,1.333)--(2.041,-3.333)--(3.674,1.333);
\draw[very thick,draw=blue](1.225,1.833)--(1.304,1.333);
\filldraw[fill=black!60,fill opacity=0.8,draw=none](2.449,1.167)--(1.633,1.333)--(1.225,1.833)--(1.633,2.167)--cycle;
\draw[draw=black,thick](1.225,1.833)--(1.633,2.167)--(2.449,1.167)--(1.633,1.333);
\filldraw[fill=black!60,fill opacity=0.8,draw=none](1.633,1.333)--(2.449,1.167)--(3.674,1.333)--cycle;
\draw[draw=black,thick](1.633,1.333)--(2.449,1.167)--(3.674,1.333);
\filldraw[fill=black!60,fill opacity=0.8,draw=none](3.674,8)--(4.082,10.833)--(6.124,8.333)--(5.307,6)--cycle;
\draw[draw=black,thick](3.674,8)--(4.082,10.833)--(6.124,8.333)--(5.307,6);
\draw[very thick](1.633,2.167)--(1.225,7.667);
\draw[very thick,draw=blue](.408,7)--(1.225,1.833);
\filldraw(.408,7) circle (2pt);
\draw[very thick,draw=blue](-.204,10.875)--(.408,7);
\filldraw[draw=black,thick,fill=black!60,fill opacity=0.8](-1.021,8.75)--(1.021,10.417)--(1.633,2.167)--(.816,1.5)--cycle;
\filldraw[fill=red!50,fill opacity=0.5,draw=none](5.307,6)--(.408,7)--(1.225,7.667)--(3.674,8)--cycle;
\draw[thick](.408,7)--(1.225,7.667)--(3.674,8)--(5.307,6);
\filldraw[fill=blue!50,fill opacity=0.5,draw=none](.408,7)--(1.304,1.333)--(3.674,1.333)--(5.307,6)--cycle;
\draw[draw=black,thin](.408,7)--(1.304,1.333);
\draw[draw=black,thin](3.674,1.333)--(5.307,6);
\filldraw[fill=black!60,fill opacity=0.8,draw=none](1.633,1.333)--(.816,1.5)--(1.225,1.833)--cycle;
\draw[draw=black,thick](1.633,1.333)--(.816,1.5)--(1.225,1.833);
\filldraw[fill=black!60,fill opacity=0.8,draw=none](1.633,1.333)--(3.674,1.333)--(.816,1.5)--cycle;
\draw[draw=black,thick](3.674,1.333)--(.816,1.5)--(1.633,1.333);
\filldraw[fill=blue!50,fill opacity=0.5,draw=none](1.633,1.333)--(3.062,-.417)--(3.674,1.333)--cycle;
\draw[draw=black,thin](3.062,-.417)--(3.674,1.333);
\filldraw[fill=black!60,fill opacity=0.8,draw=none](5.307,6)--(-.408,6.333)--(.816,1.5)--(3.674,1.333)--cycle;
\draw[draw=black,thick](-.408,6.333)--(.816,1.5)--(3.674,1.333)--(5.307,6);
\draw[very thick](.816,1.5)--(-.408,6.333);
\draw[thick](5.307,6)--(-.408,6.333);
\filldraw[fill=red!50,fill opacity=0.5,draw=none](5.307,6)--(3.674,4.667)--(-1.225,5.667)--(.408,7)--cycle;
\draw[draw=black,thin](3.674,4.667)--(-1.225,5.667);
\draw[very thick](.408,7)--(1.225,7.667);
\filldraw(1.225,7.667) circle (2pt);
\draw[very thick](1.225,7.667)--(.919,11.792);
\filldraw[fill=blue!50,fill opacity=0.5,draw=none](0,9.583)--(.408,7)--(5.307,6)--(6.124,8.333)--cycle;
\draw[draw=black,thin](5.307,6)--(6.124,8.333)--(0,9.583)--(.408,7);
\draw[very thick](-.408,6.333)--(.408,7);
\draw[thick](-.408,6.333)--(.408,7);
\filldraw[fill=black!60,fill opacity=0.8,draw=none](5.307,6)--(6.124,8.333)--(-1.021,8.75)--(-.408,6.333)--cycle;
\draw[draw=black,thick](5.307,6)--(6.124,8.333)--(-1.021,8.75)--(-.408,6.333);
\filldraw[draw=black,thin,fill=green!50,fill opacity=0.3](4.491,7)--(8.981,10.667)--(4.899,-1)--(.408,-4.667)--cycle;
\filldraw(2.041,-3.333) circle (2pt);
\filldraw[draw=white!0,fill opacity=0](-2.449,.5) circle (2pt);
\draw[very thick,draw=blue](3.674,1.333)--(5.307,6);
\filldraw[fill=red!50,fill opacity=0.5,draw=none](3.674,4.667)--(8.165,8.333)--(9.798,8)--(5.307,4.333)--cycle;
\draw[draw=black,thin](8.165,8.333)--(9.798,8)--(5.307,4.333)--(3.674,4.667);
\draw[very thick](5.307,6)--(8.165,8.333);
\draw[very thick,draw=white,dashed](0,9.583)--(1.021,10.417)--(4.082,10.833)--(6.124,8.333);
\draw[very thick,draw=white,dashed](-1.021,8.75)--(0,9.583);
\filldraw(-.408,6.333) circle (2pt);
\draw[very thick](-1.225,5.667)--(-.408,6.333);
\draw[very thick](-.408,6.333)--(-1.327,9.958);
\filldraw[draw=white!0,fill opacity=0](8.981,-1.833) circle (2pt);
\draw[very thick](3.674,4.667)--(5.307,6);
\filldraw(5.307,6) circle (2pt);
\draw[very thick,draw=blue](5.307,6)--(6.532,9.5);
\draw[very thick,draw=white,dashed](6.124,8.333)--(-1.021,8.75);
\path (4.082,10.833) node[right]{$\EuScript K$} (9.798,8) node[below right]{$\Pi$} (-1.327,9.958) node[left]{$\ell_1$} (.919,11.792) node[left]{$\ell_2$} (6.532,9.5) node[left]{$\ell$} (2.3,-3.333) node[below]{$q_0$} (2.858,11.917) node[right]{$H$} (8.981,10.667) node[right]{$H_0$} (1.9,0) node[right]{$W$} (2.449,6.167) node[below]{$\EuScript P$}(5.307,6) node[below right]{$p_0$} (-.204,10.875) node[left]{$q\vec{w}_1$}(.408,7) node[below right]{$q$} (3.4,3.6) node[above right]{$\Lambda_b$} (-1.6,5.5) node[below]{$\Lambda_a$};
\end{tikzpicture}}
&
\resizebox{7.7cm}{9cm}{
\begin{tikzpicture}[line join=round,background rectangle/.style=
{double, thick, draw=gray},
show background rectangle,xscale=0.8,yscale=0.81]
\tikzstyle{planelinestyle} = [cull=false,red,dashed]
\tikzstyle{blinestyle} = [cull=false,red,dashed]
\tikzstyle{plinestyle} = [cull=false,red,dashed]
\tikzstyle{verstyle}=[cull=false,fill=blue!20,fill opacity=0.8]
\tikzstyle{horstyle}=[cull=false,fill=blue!20,fill opacity=0.8]
\tikzstyle{boxstyle}=[cull=false,fill=blue!20,fill opacity=0.8]
\tikzstyle{polystyle} = [cull=false,fill=blue!20,fill opacity=0.8]
\filldraw[fill=green!50,fill opacity=0.5,draw=none](9.128,-4.514)--(9.128,-.629)--(1.69,.171)--(1.69,-3.714)--cycle;
\draw[draw=black,thin](1.69,.171)--(1.69,-3.714)--(9.128,-4.514);
\filldraw[fill=black!60,fill opacity=0.8,draw=none](9.128,-.629)--(4.057,-3.086)--(1.014,-1.257)--(3.55,-.029)--cycle;
\draw[draw=black,thick](4.057,-3.086)--(1.014,-1.257)--(3.55,-.029)--(9.128,-.629);
\filldraw[fill=black!60,fill opacity=0.8,draw=none](.972,-.193)--(3.55,-.029)--(1.014,-1.257)--cycle;
\draw[draw=black,thick](3.55,-.029)--(1.014,-1.257)--(.972,-.193);
\filldraw[fill=black!60,fill opacity=0.8,draw=none](.972,-.193)--(1.014,-1.257)--(1.014,-.286)--cycle;
\draw[draw=black,thick](.972,-.193)--(1.014,-1.257)--(1.014,-.286);
\filldraw[fill=red!50,fill opacity=0.5,draw=none](9.128,-4.514)--(9.128,-.629)--(9.973,-1.514)--(9.973,-4.429)--cycle;
\draw[draw=black,thin](9.973,-1.514)--(9.973,-4.429)--(9.128,-4.514);
\filldraw[fill=red!50,fill opacity=0.5,draw=none](9.128,-.629)--(9.973,-.543)--(9.973,-1.514)--cycle;
\draw[draw=black,thin](9.973,-.543)--(9.973,-1.514);
\filldraw[fill=green!50,fill opacity=0.5,draw=none](9.128,-4.514)--(11.917,-4.814)--(11.917,-.929)--(9.128,-.629)--cycle;
\draw[draw=black,thin](9.128,-4.514)--(11.917,-4.814)--(11.917,-.929);
\filldraw[fill=red!50,fill opacity=0.5,draw=none](9.128,-4.514)--(4.057,-5.029)--(4.057,-1.143)--(9.128,-.629)--cycle;
\draw[draw=black,thin](9.128,-4.514)--(4.057,-5.029);
\draw[very thick](9.128,-4.514)--(9.128,-.629);
\filldraw[fill=black!60,fill opacity=0.8,draw=none](9.128,-.629)--(11.917,-.929)--(5.578,-4)--(4.057,-3.086)--cycle;
\draw[draw=black,thick](9.128,-.629)--(11.917,-.929)--(5.578,-4)--(4.057,-3.086);
\draw[thick](9.128,-.629)--(4.057,-3.086);
\filldraw[draw=black,thin,fill=blue!50,fill opacity=0.5](0,0)--(6.761,.686)--(12.339,-.886)--(5.578,-1.571)--cycle;
\filldraw[draw=white!0,fill opacity=0](0,-7.771) circle (2pt);
\filldraw[fill=red!50,fill opacity=0.5,draw=none](9.128,-.629)--(9.128,4.229)--(9.973,4.314)--(9.973,-.543)--cycle;
\draw[draw=black,thin](9.128,4.229)--(9.973,4.314)--(9.973,-.543);
\filldraw[fill=green!50,fill opacity=0.5,draw=none](11.917,-.929)--(11.917,3.929)--(1.69,5.029)--(1.69,.171)--cycle;
\draw[draw=black,thin](11.917,-.929)--(11.917,3.929)--(1.69,5.029)--(1.69,.171);
\filldraw[fill=black!60,fill opacity=0.8,draw=none](5.747,2.914)--(9.128,-.629)--(3.55,-.029)--(1.859,1.743)--cycle;
\draw[draw=black,thick](9.128,-.629)--(3.55,-.029)--(1.859,1.743)--(5.747,2.914);
\filldraw[fill=black!60,fill opacity=0.8,draw=none](.972,-.193)--(.93,.871)--(3.55,-.029)--cycle;
\draw[draw=black,thick](.972,-.193)--(.93,.871)--(3.55,-.029);
\filldraw[draw=black,thick,fill=black!60,fill opacity=0.8](.93,.871)--(1.859,1.743)--(3.55,-.029)--cycle;
\filldraw[fill=black!60,fill opacity=0.8,draw=none](.93,.871)--(.972,-.193)--(1.014,-.286)--(1.014,.686)--cycle;
\draw[draw=black,thick](1.014,-.286)--(1.014,.686)--(.93,.871)--(.972,-.193);
\filldraw[draw=black,thick,fill=black!60,fill opacity=0.8](.93,.871)--(1.014,.686)--(1.859,1.743)--cycle;
\filldraw[fill=black!60,fill opacity=0.8,draw=none](4.057,.8)--(5.747,2.914)--(1.859,1.743)--(1.014,.686)--cycle;
\draw[draw=black,thick](5.747,2.914)--(1.859,1.743)--(1.014,.686)--(4.057,.8);
\filldraw[fill=red!50,fill opacity=0.5,draw=none](9.128,-.629)--(4.057,-1.143)--(4.057,3.714)--(9.128,4.229)--cycle;
\draw[draw=black,thin](4.057,3.714)--(9.128,4.229);
\draw[thick](4.057,-1.143)--(4.057,.8)--(5.747,2.914)--(9.128,-.629);
\draw[very thick](9.128,-.629)--(9.128,4.229);
\filldraw[fill=black!60,fill opacity=0.8,draw=none](5.747,2.914)--(7.691,3.5)--(11.917,-.929)--(9.128,-.629)--cycle;
\draw[draw=black,thick](5.747,2.914)--(7.691,3.5)--(11.917,-.929)--(9.128,-.629);
\filldraw[fill=black!60,fill opacity=0.8,draw=none](4.057,.8)--(5.578,.857)--(7.691,3.5)--(5.747,2.914)--cycle;
\draw[draw=black,thick](4.057,.8)--(5.578,.857)--(7.691,3.5)--(5.747,2.914);
\draw[very thick](4.057,-1.143)--(4.057,.8);
\filldraw[draw=black,thick,fill=black!60,fill opacity=0.8](5.578,-4)--(5.578,.857)--(1.014,.686)--(1.014,-1.257)--cycle;
\draw[thick](4.057,-3.086)--(4.057,-1.143);
\filldraw(4.057,-3.086) circle (2pt);
\draw[very thick](4.057,-3.086)--(4.057,-1.143);
\filldraw(4.057,-1.143) circle (2pt);
\draw[very thick](4.057,.8)--(4.057,3.714);
\filldraw[draw=black,thin,fill=green!50,fill opacity=0.5](0,-3.886)--(0,4.857)--(5.578,3.286)--(5.578,-5.457)--cycle;
\filldraw[draw=lightgray,fill opacity=0](0,8.743) circle (2pt);
\draw[very thick,arrows=->](4.057,-1.143)--(2.028,-.571);
\filldraw(9.128,-.629) circle (2pt);
\draw[very thick,draw=white,dashed](5.578,-1.571)--(5.578,.857)--(7.691,3.5)--(11.917,-.929);
\draw[very thick,draw=white,dashed](11.917,-.929)--(5.578,-4);
\filldraw[fill=red!50,fill opacity=0.5,draw=none](4.057,-5.029)--(3.212,-5.114)--(3.212,-1.229)--(4.057,-1.143)--cycle;
\draw[draw=black,thin](4.057,-5.029)--(3.212,-5.114)--(3.212,-1.229);
\draw[very thick](4.057,-5.029)--(4.057,-3.086);
\filldraw[fill=red!50,fill opacity=0.5,draw=none](4.057,-1.143)--(3.212,-1.229)--(3.212,3.629)--(4.057,3.714)--cycle;
\draw[draw=black,thin](3.212,-1.229)--(3.212,3.629)--(4.057,3.714);
\draw[very thick,arrows=-](4.057,.8)--(4.057,-1.143)--(2.366,-1.314);
\filldraw(4.057,.8) circle (2pt);
\draw[very thick,arrows=<-](4.057,2.743)--(4.057,.8);
\draw[very thick,arrows=->](2.366,-1.314)--(-1.014,-1.657);
\draw[very thick,draw=white,dashed](5.578,-4)--(5.578,-1.571);
\path (9.128,-.629) node[above right]{$p'_0$} (7.437,1.143) node[above right]{$\p$} (5.578,-3.709) node[above left]{$\ell'_1$} (5.578,.857) node[below left]{$\ell'_2$} (11.41,-.98) node[above left]{$\ell'$} (6.592,-4.771) node[below]{$\Pi$} (0,4.857) node[left]{$H'$} (1.69,5.029) node[above]{$H'_0$}(7.2,-1.3) node[above]{$W'$} (4.057,-1.143) node[right]{$q'$} (9.128,-4.514) node[below]{$\Lambda_b'$} (8.536,-2.243) node{$\pol'$};
\path (4.057,2.743) node[above left]{$x_3$} (-1.014,-1.657) node[below]{$x_2$} (2.028,-.571) node[above]{$x_1$};
\end{tikzpicture}}\\
(a)&(b)
\end{tabular}
\end{center}

\vspace{-3mm}
\caption{{\small First trimming position for unbounded convex polyhedra.}}\label{prirec}
\end{figure}

\vspace{1mm}
\noindent{\bf Step 5.}
Let $\pi:\R^3\to\R^3$ be the projection onto $W$ in the direction of $\vec{H}\cap\vec{\Pi}$. Then: \em $\pi(\pol\cap\Pi)=(\pol\cap\Pi)\cap W$ and $\vec{\pi}(\conv{\pol}{})=\conv{\pol}{}\cap\vec{W}$\em. Indeed, as $p_0\in\Lambda_b\cap\Pi$, $q\in\Lambda_a\cap\Pi$ and these lines are parallel to $\vec{H}\cap\vec{\Pi}$, the line $\Pi\cap W$ (through the points $p_0$ and $q$) satisfies $\pi(\pol\cap\Pi)=(\pol\cap\Pi)\cap W$. Next, we check that: $\vec{\pi}(\conv{\pol}{})=\conv{\pol}{}\cap\vec{W}$.
 
Notice first that $\vec{H}^+\cap\vec{H}_0^+\cap\vec{\Pi}^+\cap\vec{W}$ is the cone generated by two vectors $\vec{w}_1,\vec{w}_2\in\conv{\pol}{}\setminus\{0\}$ (non-necessarily linearly indepedent) such that $\vec{w}_1$ generates the line $W\cap H$ and $\vec{w}_2$ the line $H_0\cap W$. Thus,\setcounter{equation}{0}
\begin{equation}\label{1} 
\vec{H}^+\cap\vec{H}_0^+\cap\vec{\Pi}^+\cap\vec{W}\subset\conv{\pol}{}\cap\vec{W}\subset\vec{\pi}(\conv{\pol}{}).
\end{equation}
On the other hand, since $\pol$ is in FU-position w.r.t. $\vec{\Pi}$ and $\pol\subset H^+\cap H_0^+$, we deduce that
$$
\conv{\pol}{}\subset\conv{H^+\cap H_0^+}{}\cap\vec{\Pi}^+=\vec{H}^+\cap\vec{H}_0^+\cap\vec{\Pi}^+.
$$
As $\vec{H}\cap\vec{\Pi}\subset\vec{H}_0$ and $\im(\vec{\pi})=\vec{W}$, we have
\begin{equation}\label{2}
\vec{\pi}(\conv{\pol}{})\subset\vec{H}^+\cap\vec{H}_0^+\cap\vec{\Pi}^+\cap\vec{W}.
\end{equation}
Combinining \eqref{1} and \eqref{2}, we conclude that $\vec{\pi}(\conv{\pol}{})=\conv{\pol}{}\cap\vec{W}$.

\vspace{1mm}
\noindent{\bf Step 6.}
After a change of coordinates that transforms $H^+$ in $\{x_2\leq0\}$, $W$ in $\{x_3=0\}$ and the line $\vec{H}\cap\vec{\Pi}$ onto the line $\{x_1=x_2=0\}$, we conclude by Lemma \ref{proysec} and Remark \ref{sectioning} that the convex polyhedron $\pol$ is in first trimming position with respect to the facet $\Ff$, as required.
\end{proof}

\begin{proof}[Proof of Proposition \em \ref{1trim} (ii) \em for bounded facets]
We check first that $\pol$ has a bounded facet. Assume that the unbounded edges of $\pol$ are parallel to the vector $\vec{e}_3$. Let $\Aa:=p\vec{e}_3$ be an unbounded edge of $\pol$. Let $\Ff$ be a facet of $\pol$ non-parallel to $\vec{e}_3$ such that $p\in\Ff$. It holds that $\conv{\Ff}{}\subset\conv{\pol}{}=\{\lambda\vec{e}_3:\ \lambda\geq0\}$. As $\Ff$ is non-parallel to $\vec{e}_3$, we have $\conv{\Ff}{}=\{{\bf0}\}$, so $\Ff$ is bounded.

Fix a bounded facet $\Ff_0$ of $\pol$ and let $H\subset\R^3$ be the plane generated by $\Ff_0$. Since the unbounded edges of $\pol$ are parallel, $\conv{\pol}{}=\{\lambda\vec{v}:\ \lambda\geq0\}$. Let $\vec{\ell}$ be the line generated by $\vec{v}$. As $\Ff_0=\pol\cap H$ is a bounded facet, $\vec{H}\cap\vec{\ell}=\{{\bf0}\}$. Let now $\Pi$ be a plane parallel to $H$ that meets all the unbounded edges of $\pol$ and such that all the bounded faces of $\pol$ are contained in $\Int\Pi^-$. By Remark \ref{sectioning} $\p:=\pol\cap\Pi$ is a bounded convex polygon, $\pol\cap\Pi^-$ is a bounded convex polyhedron and $\pol\cap\Pi^+=\p+\conv{\pol}{}=\p+\{\lambda\vec{v}:\ \lambda\geq0\}$. Thus, $\pol\cap\Pi^+$ is affinely equivalent to $\p\times{[0,{+\infty}[}$. After a change of coordinates, we assume that:
\begin{itemize}
\item $H$ is the plane $\{x_2=0\}$ and $\pol\subset\{x_2\leq0\}$,
\item $\Pi$ is the plane $\{x_2=-1\}$ and $\pi_3(\p)=\p\cap\{x_3=0\}$ (recall here the projection property for convex polygons described in the Introduction),
\item $\vec{v}=-\vec{e}_2$, so $\conv{\pol}{}=\{-\lambda\vec{e}_2:\ \lambda\geq0\}\subset\{x_3=0\}$ and $\pi_3(\conv{\pol}{})=\conv{\pol}{}\cap\{x_3=0\}$. 
\end{itemize}
By Lemma \ref{proysec} and Remark \ref{sectioning} $\pol$ is in first trimming position with respect to the facet $\Ff_0$.
\end{proof}

\subsection{Second trimming position}

Our purpose now is to prove Proposition \ref{2trim}.

\begin{proof}[Proof of Proposition \em \ref{2trim}]
(i) Assume that $\pol\cap\{x_3=0\}=\Ff$ is one of its unbounded facets with non-parallel unbounded edges $\Aa_1$ and $\Aa_2$. Let $H_i$ be the plane of $\R^3$ generated by the facet $\Ff_i$ such that $\Aa_i=\Ff\cap\Ff_i$. As $H_1\cap\{x_3=0\}$ and $H_2\cap\{x_3=0\}$ are non-parallel lines, we may assume keeping invariant the plane $\{x_3=0\}$ that $H_i:=\{x_i=0\}$. Changing the signs of the variables if necessary we assume that $\pol\subset\{x_1\geq0,x_2\geq0,x_3\leq0\}$. It remains to show that ${\mathfrak A}_\pol$ is bounded. Consider now the projection $\pi_3:\R^3\to\R^2, (x_1,x_2,x_3)\mapsto (x_1,x_2,0)$.

It is clear that $\mathfrak A_\pol\subset \pi_3(\pol)\subset\{x_1\ge0, x_2\ge0\}$. Take the extremes $p_1=(0, b, 0)$, $p_2=(a,0,0)$ of the unbounded edges $\Aa_1$ and $\Aa_2$, where $a,b\ge 0$. The bounded convex polygon 
$$
\left\{(x_1,x_2)\in\R^2: x_1\ge 0, x_2\ge 0, \frac{x_1}{a+1}+\frac{x_2}{b+1}-1\le 0\right\}
$$
contains $\mathfrak A_\pol$, so it is bounded. Thus, $\pol$ is in second trimming position with respect to $\Ff$.

(ii) If $\Ff$ is a bounded facet of $\pol$, then $\Ff$ is non-parallel to the unbounded edges of $\pol$. We may assume that $\pol\subset\{x_3\leq 0\}$, $\pol\cap\{x_3=0\}=\Ff$ is a facet of $\pol$ and the unbounded edges of $\pol$ are parallel to the vector $\vec{e}_3$. Thus, $\pol$ is in FU-position w.r.t. $\vec{\Pi}:=\{-x_3=0\}$. By \ref{conoedges} and as $\vec{\pi}_3(\vec{e}_3)={\bf0}$, we deduce that $\pi_3(\pol)$ is bounded, so $\mathfrak A\subset \pi_3(\pol)$ is also bounded. Thus, $\pol$ is in second trimming position with respect to $\Ff$.
\end{proof}

\appendix
\section{Limitations of the trimming positions}\label{s4}

In this appendix we construct unbounded convex polyhedra that cannot be placed in first or second trimming position with respect to their facets. The clue is that the Property stated at the end of the Introduction fails for $n>2$. The key result is the following, which can be proved by perturbing appropriately the vertices of a \em regular cross-polytope \em of $\R^n$, which is the convex hull of the set of points consisting of all the permutations of $(\pm1,0,\ldots,0)\in\R^n$.

\begin{lem}[Convex polyhedra without sectional projections]\label{nsp}
For each $n\geq3$ there exists an $n$-dimensional bounded convex polyhedron $\pol\subset\R^n$ such that for each line $\vec{\ell}$ and each hyperplane $H$ non-parallel to $\vec{\ell}$ it holds $\pi(\pol)\neq\pol\cap H$ where $\pi:\R^n\to H$ is the linear projection onto $H$ in the direction of $\vec{\ell}$.
\end{lem}

\begin{cor}[Counterexamples to the trimming positions]\label{chungoppr}
Let $n\geq4$ and let $\p\subset\R^{n-1}$ be a bounded convex polyhedron without sectional projections. Define $\pol:=\p\times{[0,{+\infty}[}\subset\R^n$. Then $\pol$ can be placed neither in first trimming position with respect to any of its facets nor in second trimming position with respect to any of its unbounded facets.
\end{cor}
\begin{proof}
Denote $\Ff_0:=\pol\cap\{x_n=0\}=\p\times\{0\}$, which is a facet of $\pol$ and let $\vec{\ell}$ be the direction of the unbounded edges of $\pol$. Suppose that $\pol$ is placed in first trimming position with respect to one of its facets. By Lemma \ref{proysec} there exists a hyperplane $\Pi\subset\R^n$ such that: its direction $\vec{\Pi}$ contains the line $\vec{\ell}_n$ generated by $\vec{e}_n$, $\Pi$ meets each unbounded edge of $\pol$ in a singleton, the vertices of $\pol$ are contained in $\Int{\Pi}^-$ and $\pi_n(\pol\cap\Pi)=\pol\cap\Pi\cap\{x_n=0\}$. \setcounter{substep}{0}

\begin{substeps}{chungoppr}\label{finalprueba}
Consider the projection $\rho:\R^n\to\R^n$ in the direction of the line $\vec{\ell}$ onto the hyperplane $H$ generated by $\Ff_0$. As $\Pi$ meets the unbounded edges of $\pol$ in a singleton and all of them are parallel to $\vec{\ell}$, we deduce that $\Pi$ is non-parallel to $\vec{\ell}$, so $g:=\rho|_{\Pi}:\Pi\to H$ is an affine bijection. By Remark \ref{sectioning} $\p_1:=\pol\cap\Pi$ is the bounded convex polyhedron whose vertices are the intersections of the unbounded edges of $\pol$ (all of them parallel to $\vec{\ell}$) with the hyperplane $\Pi$. As $\vec{\rho}(\vec{\ell})=\{{\bf0}\}$, we deduce that $\rho(\p_1)=\Ff_0$. If we define $\vec{r}:=\vec{g}(\vec{\ell}_n)$ and $W:=g(\Pi\cap\{x_n=0\})\subset H$, the projection $\pi:=(g\circ\pi_n\circ g^{-1}):H\to H$ in the direction of $\vec{r}$ onto the hyperplane $W$ satisfies $\pi(\Ff_0)=\Ff_0\cap W$, against the fact that $\Ff_0:=\p\times\{0\}$ has no sectional projections. Thus, $\pol$ cannot be placed in first trimming position with respect to any of its facets.
\end{substeps}

Suppose next that $\pol$ is placed in second trimming position with respect to one of its unbounded facets $\Ff_0$, which is contained in the hyperplane $\{x_n=0\}$. The unbounded edges of $\pol$ are parallel to a line $\vec{\ell}$ contained in the hyperplane $\{x_n=0\}$, so they are non-parallel to the line $\vec{\ell}_n$ generated by $\vec{e}_n$. By Lemma \ref{proysec} there exists a hyperplane $\Pi\subset\R^n$ such that: its direction $\vec{\Pi}$ contains the line $\vec{\ell}_n$, the hyperplane $\Pi$ meets the unbounded edges of $\pol$ (which are all non-parallel to $\vec{\ell}_n$), the bounded edges of $\pol$ are contained in the open half-space $\Int{\Pi}^-$ and $\pi_n(\pol\cap\Pi)=\pol\cap\Pi\cap\{x_n=0\}$. Proceeding analogously to (\ref{chungoppr}.\ref{finalprueba}) we achieve a contradiction. Thus, $\pol$ cannot be placed in second trimming position with respect to any of its unbounded facets.
\end{proof}

\end{document}